%% file: itsc_final.tex
\documentclass[final,conference]{IEEEtran}
\IEEEoverridecommandlockouts
\usepackage{cite}
\usepackage{amsmath,amssymb,amsfonts}
\usepackage{algorithmic}
\usepackage{graphicx}
\usepackage{textcomp}
\usepackage{xcolor}
\usepackage{cleveref}
\usepackage{glossaries}
\usepackage{pgfplots}
\usepackage{pgfplotstable}
\usepackage{xurl}   

\pgfplotsset{compat=1.18}

\input{tex/param-barplots.tex}
\newacronym{iam}{IAM}{Innovative Air Mobility}
\newacronym{blq}{BLQ}{Bologna Guglielmo Marconi Airport}
\def\BibTeX{{\rm B\kern-.05em{\sc i\kern-.025em b}\kern-.08em
    T\kern-.1667em\lower.7ex\hbox{E}\kern-.125emX}}

\begin{document}

\title{Short‑Range IAM Operations to Overcome Access Gaps Around Bologna Airport}

\author{
    \IEEEauthorblockN{1\textsuperscript{st} Jean-Claude Lebègue}
    \IEEEauthorblockA{\textit{OPTIM} \\
    \textit{Ecole Nationale de l'Aviation Civile}\\
    Toulouse, France \\
    jean-claude.lebegue@enac.fr}
    \and
    \IEEEauthorblockN{2\textsuperscript{nd} Xuhao Gui}
    \IEEEauthorblockA{\textit{OPTIM} \\
    \textit{Ecole Nationale de l'Aviation Civile}\\
    Toulouse, France \\
    xuhao.gui@enac.fr}
    \and
    \IEEEauthorblockN{3\textsuperscript{rd} Daniel Delahaye}
    \IEEEauthorblockA{\textit{OPTIM} \\
    \textit{Ecole Nationale de l'Aviation Civile}\\
    Toulouse, France \\
    daniel.delahaye@enac.fr}
}

\maketitle

\begin{abstract}
We study \gls{iam} as a dedicated hub feeder and dispersal service for \gls{blq}, addressing short‑range airport‑to‑airport connectivity where scheduled fixed‑wing services are inefficient and surface access is time‑consuming. Demand is generated directly from \gls{blq}’s flight timetable via a staggered rolling‑window model: arriving passengers become \gls{iam}‑eligible after a 20‑minute buffer (a baseline 70\% “staying” share is considered and distributed across nearby airports using a gravity model), while feeder demand for departures is timestamped backward by \gls{iam} travel time plus a 2‑hour pre‑flight buffer. We build an \gls{iam}‑only time‑expanded network (airports as vertiports) and solve two sequential subproblems: passenger routing with A* to reveal edge‑level demand pressure, followed by fleet dispatch with dynamic programming to schedule a finite \gls{iam} fleet. Using OAG schedules for 5 May 2024, we simulate 09:00–17:00 with a 50‑vehicle fleet (10 seats, 30 m/s) and 5‑minute windows. Results show large travel‑time advantages versus rail for many links (e.g., \gls{blq}–Il Ciocco: 46 min \gls{iam} vs. 237–294 min by train), while a few pairs remain competitive for rail (e.g., \gls{blq}–Florence). Sensitivity to \gls{iam} speed and staying share quantifies how technology and adoption shape feasible coverage, timing feasibility, and hub accessibility, demonstrating the potential of \gls{iam} to strengthen \gls{blq}’s regional reach without altering airline schedules.
\end{abstract}


\section{Introduction}

Over the past decade, European transport policy has emphasized improving door‑to‑door mobility, notably through Flightpath2050~\cite{flightpath2050}, which targets four‑hour door‑to‑door journeys for most travelers. Achieving this vision requires addressing access gaps around major hubs where nearby airports are too close to justify scheduled fixed‑wing flights, yet far enough that surface access remains time‑consuming. In such settings, further expansion of heavy ground infrastructure can be economically and geographically challenging.

\gls{iam}—using new generations of vertical or short‑takeoff aircraft—offers a targeted alternative~\cite{easa2021study} for short ranges, enabling fast, flexible connections with limited infrastructure. Prior European initiatives have explored how \gls{iam} might complement established aviation, and our earlier work proposed both heuristic and exact optimization frameworks on time‑expanded networks to study regional \gls{iam} planning. Those studies focused on multimodal integration or general regional scenarios and highlighted both the potential of \gls{iam} and the computational challenges of large‑scale planning.

This paper shifts the scope to airport‑to‑airport connectivity with \gls{iam} as the only modeled transport mode. We consider \gls{blq} and evaluate \gls{iam} as a hub feeder/dispersal system linking \gls{blq} to nearby intra‑ and inter‑regional airports. In our framework, \gls{blq} is the only node with a commercial flight schedule, which is used solely to time‑stamp demand (arrivals for dispersal and departures for feeders). The network itself is an \gls{iam}‑only time‑expanded graph in which nearby airports are treated as \gls{iam}‑capable facilities (vertiports). The goal is to quantify how \gls{iam} can improve \gls{blq}’s regional accessibility without altering airline schedules.

Methodologically, we use a two‑stage, sequentially coupled approach. First, passenger routing applies A* on the \gls{iam}‑only time‑expanded network under an unconstrained \gls{iam} assumption to reveal demand‑pressure on \gls{iam} legs. Second, dynamic programming schedules a finite \gls{iam} fleet over the same graph using those pressures as edge rewards. Demand is generated from \gls{blq}’s schedule via a staggered rolling window: post‑arrival dispersal after a 20‑minute buffer, of which a share (passenger diversion) uses \gls{iam} and is distributed to \gls{iam} airports via a gravity‑model based on travel time and population; and feeder demand for departures, time‑stamped backward by \gls{iam} travel time plus a pre‑flight buffer. Using real flight data, we then assess the regional value of \gls{iam} and examine how performance varies with two key factors—\gls{iam} cruise speed and the passenger diversion share.

\Cref{sec:related} reviews related work with emphasis on the shift from multimodal \gls{iam} planning to \gls{iam}‑only hub feeders. \Cref{sec:methodology} presents the methodology: framework overview, demand collection, \gls{iam}‑only time‑expanded network, passenger routing with A*, and dispatch via dynamic programming. \Cref{sec:results} reports results for the \gls{blq} case. \Cref{sec:conclusion} concludes.

\section{Related work} \label{sec:related}

Research on \gls{iam} has expanded rapidly in recent years~\cite{li2025evaluation,wu2025optimization}, driven by increasing interest in integrating new-generation vertical or short‑takeoff aircraft into multimodal transport systems. Early studies have focused primarily on airspace organization, urban operations, congestion management, or vertiport placement models within the broader Urban Air Mobility (UAM) context~\cite{chati2025vertiport}. While these works provide foundational insights for emerging air mobility ecosystems, they concentrate mainly on intra‑urban settings and do not directly address the challenges of short‑range airport‑to‑airport connectivity in regional environments.
A complementary line of research has begun to analyze \gls{iam} as a means of enhancing regional accessibility, especially in rural or geographically constrained areas. Our earlier contributions are part of this evolving literature. In our first study, we proposed a regional‑scale multimodal \gls{iam} service optimization framework built on a heuristic architecture combining A* search and Dynamic Programming (DP)~\cite{gui2025regional}. The key idea was to decompose the problem into two interdependent NP‑hard subproblems:

\begin{itemize}
    \item Passenger journey optimization, formulated as a multi‑commodity flow problem  on a time‑expanded transport network~\cite{GUITART2026103363}, where A* search was used to compute feasible multimodal passenger paths and estimate edge-level demand pressures;
    \item \gls{iam} fleet dispatching, formulated as a trajectory‑planning problem analogous to a Directed Steiner Tree variant~\cite{charikar1999approximation}, where DP was used to generate \gls{iam} flight schedules that respond to the demand structure identified in the passenger‑routing phase.
\end{itemize}

This separation allowed the framework to remain computationally tractable while capturing the interactions between passenger behavior and \gls{iam} fleet operations.

In a subsequent study, we extended this methodological foundation by formulating the \gls{iam} dispatching problem as a Mixed‑Integer Linear Programming (MILP) model~\cite{gui:hal-05531747}. This exact formulation explicitly integrated passenger routing and \gls{iam} fleet movements in a single optimization structure, providing globally optimal solutions for small‑ to medium‑scale networks. To mitigate the combinatorial explosion inherent in MILP models, we introduced a vertiport‑connection sparsification technique, drastically reducing the number of air‑transport edges while preserving essential connectivity patterns. Although the MILP model offered mathematical rigor and served as a valuable benchmark, its scalability limitations highlighted the need for lighter, more flexible optimization methods when large networks or near‑real‑time decision support are required.

Together, these two previous studies established a methodological framework for regional‑scale \gls{iam} analysis, combining heuristic tractability with exact optimization insights. However, both works focused on general multimodal transport scenarios, without explicitly examining the structural constraints of closely spaced airports—a setting where traditional commercial flights are inefficient but ground travel remains slow. This specific operational gap has received little attention in the \gls{iam} literature, and the potential role of \gls{iam} as a feeder mode within dense regional air networks, particularly for enhancing the reach of major hubs, remains largely unexplored.
The present paper extends this line of research by applying the two‑subproblem modeling paradigm to a new and operationally relevant context: airport‑to‑airport connectivity in northern Italy, with a focus on strengthening feeder links to \gls{blq}. By adapting the time‑expanded network modeling approach and the A*+DP heuristic framework to this new use case, the study provides the first dedicated assessment of \gls{iam} as a replacement for inefficient short‑haul flights, filling an important gap in the current literature and expanding the applicability of the \gls{iam} optimization methodologies developed in our earlier works.

\section{Methodology} \label{sec:methodology}

This section presents the methodology to evaluate \gls{iam} as a short‑range hub feeder/dispersal service around \gls{blq} using an \gls{iam}‑only time‑expanded network. The problem is structured into two sequential subproblems: passenger routing on the \gls{iam}‑only graph to reveal time‑dependent demand pressure on \gls{iam} links, and \gls{iam} dispatching to schedule a finite fleet against those passenger‑driven priorities. We first outline the \gls{blq}‑centric framework and demand model, then describe the network construction and both optimization stages.

\subsection{Framework overview}

The methodology adopts a two‑stage, sequentially coupled framework that remains computationally tractable while preserving a tight alignment between passenger preferences and \gls{iam} operational feasibility. First, a passenger journey optimization stage identifies time‑optimal itineraries on a time‑expanded airport network under an unconstrained \gls{iam} assumption. This assumption is intentional: it reveals the \gls{iam} legs that passengers would most value if capacity were not a bottleneck, producing edge‑level popularity weights (demand‑pressure indicators) that reflect when and where \gls{iam} links would be most beneficial. Second, an \gls{iam} dispatching optimization stage uses those popularity weights as edge rewards to schedule a finite \gls{iam} fleet on the same time graph, assigning aircraft to high‑value feeder legs while enforcing vehicle capacity, duty‑time, and temporal consistency constraints. After each vehicle is scheduled, the framework reduces the residual demand on traversed legs and iterates to the next vehicle, ensuring that subsequent assignments continue to target unmet, high‑priority demand. In this way, passenger routing informs fleet decisions and fleet decisions feed back into remaining demand, delivering demand‑consistent schedules without the computational burden of a fully joint (and intractable) optimization. The design also lends itself to scenario analysis—varying fleet size, \gls{iam} performance, or airport sets—because the passenger stage can be recomputed quickly to update demand signals, and the dispatching stage can be rerun to produce new schedules aligned with the latest passenger‑centric preferences.

\subsection{Demand collection and estimation}

Passenger demand in this study is generated directly from the flight schedule of the hub airport. 
Unlike our previous work---where multimodal regional demand was inferred from a combination of 
population data and public transport schedules---the present study focuses exclusively on 
airport--to--airport connectivity and therefore models demand solely from the temporal flow of 
passengers arriving at or departing from the hub. To ensure temporal consistency with \gls{iam} 
operations, the demand collection mechanism is implemented through a staggered rolling-window
strategy, inspired by the dynamic demand--supply estimation scheme used in our earlier frameworks \cite{gui2025regional, gui:hal-05531747}.

For each rolling window of fixed duration, we consider two categories of passengers. The first 
category corresponds to travelers arriving at the hub who may subsequently continue their 
journey toward nearby airports using \gls{iam}. To capture their availability for \gls{iam} transfers, we assume that passengers require approximately 20 minutes to disembark, retrieve their belongings, and reach the \gls{iam} boarding area. Therefore, for each window, we examine the hub's arrival schedule 20 minutes earlier and identify the passengers whose flights landed in that preceding window. For example, a flight arriving at the hub at time $t_{arr}$ generates \gls{iam} demand at time $t$, defined as:
\begin{equation}
t = t_{\mathrm{arr}} + \tau_{\mathrm{dbk}},
\end{equation}
where $\tau_{dbk}$ denotes the disembarkation time. The resulting demand is then assigned to the time window containing $t$.

Only a fraction $p \in [0,1]$ of the passengers arriving in a time window is assumed to transfer via \gls{iam}, with the remaining proportion connecting to other commercial flights at the hub. To distribute this \gls{iam}-attributed demand across the set of nearby \gls{iam}-enabled airports, a gravity-based model is employed. Each airport is associated with a population catchment area and an \gls{iam} travel time; airports serving larger populations and reachable more quickly receive a proportionally higher share of the dispersal demand. 


Let us consider an arrival time window $\omega_{arr}$. Let $d_{\omega_{arr}}$ denote the number of passengers arriving at the hub $H$ during this window. A fraction $p$ of these passengers is assumed to generate \gls{iam} dispersal demand. The total \gls{iam} demand $O_\omega$ assigned to the \gls{iam} network for window $\omega$ is therefore defined as:
\begin{equation}
    O_\omega = p \cdot d_{\omega_{\mathrm{arr}}}.
\end{equation}
This demand is then distributed from the hub $H$ to the set of nearby airports $\{A_i\}_i$ using a gravity-based model. The share of demand $T_{H,A_i}$ assigned to airport $A_i$ is given by:
\begin{equation}
    T_{H,A_i} = O_\omega
    \frac{(P_i^{\lambda} / \tau_{H,A_i}^{\mu})^{\alpha}}
    {\sum_k (P_k^{\lambda} / \tau_{H,A_k}^{\mu})^{\alpha}},
\end{equation}
where $P_i$ denotes the population of the city associated with airport $A_i$, and $\tau_{H,A_i}$ is the \gls{iam} travel time between the hub $H$ and airport $A_i$. The parameters $(\lambda, \mu, \alpha)$ control the relative influence of population size, travel time, and demand concentration in the gravity model.

The second demand category corresponds to passengers who must reach the hub $H$ in the future to 
catch an outgoing flight. For each scheduled departure, we determine the latest feasible time at 
which a passenger may leave an \gls{iam}-enabled airport $A_i$ while still arriving at the hub $H$ with sufficient time to check in and complete international procedures. Let $\tau_{A_i,H}$ denote the \gls{iam} travel time between $A_i$ and $H$, and $b$ the pre-departure buffer time. A departure scheduled at time $t_{\mathrm{dep}}$ generates \gls{iam} feeder demand at time $t$ as follows:

\begin{equation}
t = t_{\mathrm{dep}} - \tau_{A_i,H} - b.
\end{equation}

This demand is assigned to the rolling window containing $t$. This forward-looking time stamping ensures that \gls{iam} feeder movements are temporally aligned with the actual flight departure requirements.

The rolling-window mechanism consolidates these two demand streams---inbound dispersal and outbound feeder demand---into a unified and chronologically consistent sequence of \gls{iam} service requests. Although the time-expanded network models only \gls{iam} movements, the demand estimation process ensures that all \gls{iam} flows are synchronized with the temporal structure of hub flight operations. This approach extends the demand-collection logic of our previous studies while adapting it to a purely airport--to--airport setting in which the hub flight schedule defines the timing of demand, and \gls{iam} represents the sole transport mode connecting the hub to the surrounding region.

\subsection{Time‑Expanded Network Construction}

The network follows the time‑expanded representation from our earlier \gls{iam} planning work \cite{gui2025regional}, adapted here to an \gls{iam}‑only airport‑to‑airport context. The static topology consists of \gls{blq} (the hub) and a set of nearby airports modeled as \gls{iam}‑capable facilities; \gls{iam} can operate bidirectionally between \gls{blq} and each facility, forming a star‑shaped graph. Time expansion replicates each airport across discrete time layers; directed arcs encode \gls{iam} flight legs (with travel times derived from inter‑airport distance and assumed cruise speed) and waiting arcs. Because arcs connect only forward in time, the graph is a directed acyclic graph, enabling efficient A* search and dynamic programming. No commercial flight arcs are included: \gls{blq}’s flight schedule is used only to time‑stamp demand.

\subsection{Passenger journey optimization}

Passenger journey optimization determines which \gls{iam} legs passengers would prefer to use to reach or leave the hub at the required times. The search space is the \gls{iam}‑only time‑expanded airport graph comprising \gls{iam} flight arcs and waiting arcs. For each demand item generated by the rolling‑window model, the A* algorithm computes the minimum‑time path subject to temporal feasibility. Because the graph is a DAG ordered by time, A* converges quickly using distance‑based heuristics. The output is a set of time‑optimal \gls{iam} itineraries and edge‑level usage counts, which act as demand‑pressure indicators for the subsequent dispatching stage. By assuming unlimited \gls{iam} capacity at this stage, the framework isolates the intrinsic passenger preferences independently of fleet limitations.

\subsection{\gls{iam} Dispatching Optimization}

The \gls{iam} dispatching stage translates passenger‑driven demand indicators into actual \gls{iam} flight schedules that can be feasibly operated by a finite fleet. The dispatching problem is defined on the same time‑expanded network but restricted to nodes and edges relevant to \gls{iam} operations. Each \gls{iam} vehicle is modeled as an agent that can occupy exactly one airport at any given time and can traverse \gls{iam} edges only if the corresponding flight time and connectivity constraints are satisfied. Dynamic Programming (DP) is used to compute the highest‑reward feasible itinerary for each \gls{iam} vehicle. Since the time-expanded graph is acyclic, DP proceeds chronologically, accumulating rewards associated with edges that are heavily demanded according to the passenger optimization stage. After a route is assigned to an \gls{iam} vehicle, the passenger demand associated with the traversed \gls{iam} edges is reduced, ensuring that subsequent vehicles are allocated to remaining high‑value portions of the network. The process is repeated until all \gls{iam} vehicles have been dispatched. This sequential dispatching strategy strikes a balance between computational tractability and operational realism, enabling the construction of efficient \gls{iam} schedules that respond to demand variability across the regional airport system.

\section{Results} \label{sec:results}

The following section presents the results obtained by applying the proposed \gls{iam}‑only feeder framework (\cref{sec:methodology}) to the regional airport system surrounding \gls{blq}. \gls{blq} provides an ideal context for evaluating short‑range \gls{iam} connectivity because it is the principal air transport hub of the Emilia‑Romagna region and concentrates the overwhelming majority of the area’s commercial flight operations. Several nearby airports lie within distances too short to justify scheduled fixed‑wing flights, yet far enough that surface access remains time‑consuming, creating structural accessibility gaps for passengers attempting to reach or disperse from the hub. Moreover, \gls{blq} exhibits pronounced peaks in arrivals and departures throughout the day, making it a realistic testbed for synchronizing \gls{iam} operations with hub‑driven passenger flows.

\begin{table}
\centering
\caption{Comparison of travel times (in minutes) between Bologna and nearby airports/heliports by train and \gls{iam} (30 m/s).}
\label{tab:train_vs_IAM}
\resizebox{.95\columnwidth}{!}{%

\begin{tabular}{lcccc}
\hline
\textbf{ID} & \textbf{Name} & \textbf{Train (min)} & \textbf{Train (max)} & \textbf{\gls{iam} (min)} \\
\hline
K32860 & Dominio di Bagnoli Airfield              & 140 & 192 & 47 \\
K43517 & Il Ciocco Heliport                       & 237 & 294 & 46 \\
K43547 & C\`a Negra Private Airstrip              & 145 & 196 & 51 \\
K43548 & Reggio Emilia Airport                    & 37  & 49  & 29 \\
K43550 & Lugo Francesco Baracca Airfield          & 38  & 41  & 26 \\
K43554 & Legnago Airfield                         & 100 & 157 & 37 \\
K43555 & Pavullo Airfield                         & 137 & 171 & 24 \\
K43556 & Ravenna Airport                          & 72  & 106 & 42 \\
K43558 & Carpi Budrione Airfield                  & 111 & 178 & 26 \\
K43559 & Prati Vecchi d'Aguscello Airfield        & 26  & 58  & 23 \\
K43597 & Parma Airport                            & 55  & 68  & 47 \\
K43606 & Cervia Air Base                          & 37  & 74  & 49 \\
K43609 & Ferrara San Luca Airfield                & 26  & 58  & 22 \\
K43613 & Forl\`i Luigi Ridolfi Airport            & 36  & 74  & 40 \\
K43628 & Lucca Tassignano Airport                 & 153 & 175 & 54 \\
K43644 & Florence Airport (Peretola)              & 37  & 38  & 45 \\
\hline
\end{tabular}
}
\end{table}

A comparison of \gls{iam} travel times with existing train services clearly highlights the potential of \gls{iam} to enhance regional accessibility to and from \gls{blq}. \Cref{tab:train_vs_IAM} shows that \gls{iam} consistently outperforms rail travel across all considered airports and heliports, often by a substantial margin. For several locations—such as Pavullo Airfield, Carpi Budrione Airfield, Ferrara San Luca Airfield, and Prati Vecchi d'Aguscello Airfield—the \gls{iam} travel time is less than half of the fastest available train connection. Even for larger regional airports such as Parma, Ravenna, and Forlì, \gls{iam} offers competitive or significantly faster access. The difference is especially pronounced for more remote airports like Il Ciocco Heliport or Lucca Tassignano Airport, where \gls{iam} reduces travel time by more than two hours relative to conventional ground transport. 

The empirical analysis relies on real flight schedule data obtained from OAG~\cite{oag2024}, using all arrivals and departures recorded for 5 May 2024 at \gls{blq}. The simulation spans 09:00–17:00 with 5‑minute rolling windows. The \gls{iam} fleet comprises 50 vehicles with 10 seats each and a nominal cruise speed of 30 m/s. A load factor of 0.8 is applied to flights to approximate passenger volumes. A passenger diversion of 0.7 represents the proportion of passenger volume using the \gls{iam}. For outbound feeders, we impose a 2‑hour pre‑departure buffer to account for check‑in, security, and boarding at \gls{blq}. These settings define the baseline environment in which \gls{iam} demand is generated and dispatch is evaluated.

\gls{blq} typically handles around 200 flight movements per day according to publicly available statistics~\cite{airnav2024}, the dataset used in this study contains only about 50 daily flights, representing roughly one quarter of the airport’s actual activity. This discrepancy arises from limitations in the available timetable data, which does not fully capture the complete set of arrivals and departures operating at \gls{blq}. To ensure that the demand estimation reflects realistic traffic levels, we introduce a scaling factor equal to 4. This factor proportionally increases the number of arriving and departing passengers in each time window, thereby restoring the total daily flight volume to a level consistent with the known operational intensity of the airport.

\begin{figure}
    \centering
    \includegraphics[width=\linewidth]{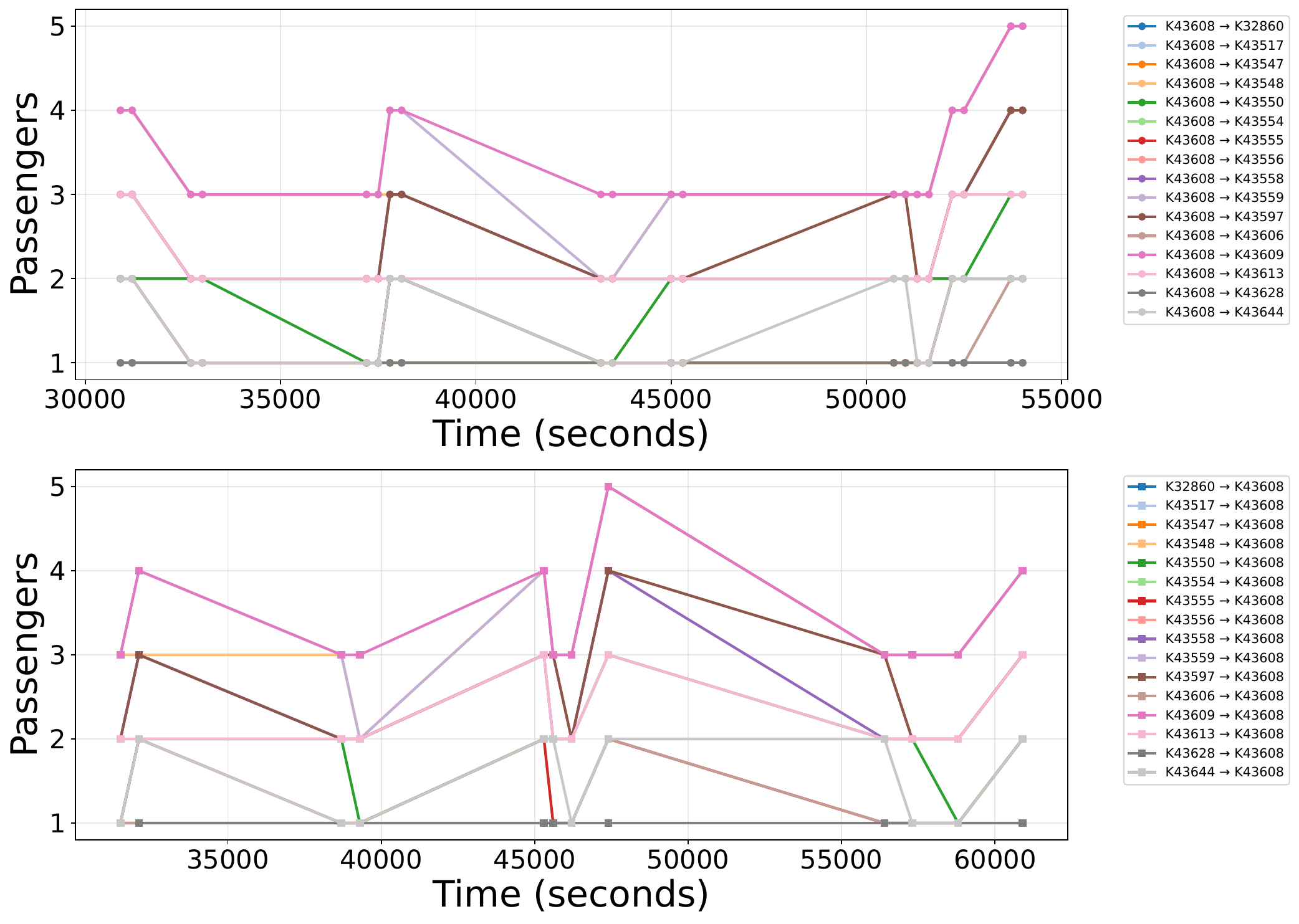}
    \caption{Simulated Time‑Series of Departing (top) and Arriving (bottom) Passenger Demand at Node K43608 (\gls{blq})}
    \label{fig:pax_timeserie}
\end{figure}

The simulated time‑series of passenger demand associated with node K43608 (\gls{blq}) reveals distinct temporal patterns for both departing and arriving flows (see \cref{fig:pax_timeserie}). Outbound demand shows moderate variability across the main OD pairs, with several relations—such as K43608→K43609, K43608→K43547, and K43608→K43548—consistently exhibiting higher passenger volumes throughout the simulated period. The demand profiles follow smooth fluctuations rather than abrupt peaks, indicating a relatively stable evolution of traffic over time. In contrast, inbound demand toward K43608 demonstrates a slightly broader distribution of values, with OD pairs like K43609→K43608 and K43547→K43608 contributing significantly to the total incoming flow. The temporal progression for arrivals also displays coherent variations, with no extreme surges, suggesting that the node operates under balanced and predictable demand dynamics. Overall, the results highlight K43608 as a well‑connected node with steady bidirectional passenger movement, where both outgoing and incoming demands evolve in a coordinated and stable manner over the simulated timeframe.

\begin{figure}
  \centering
  \includegraphics[width=\linewidth]{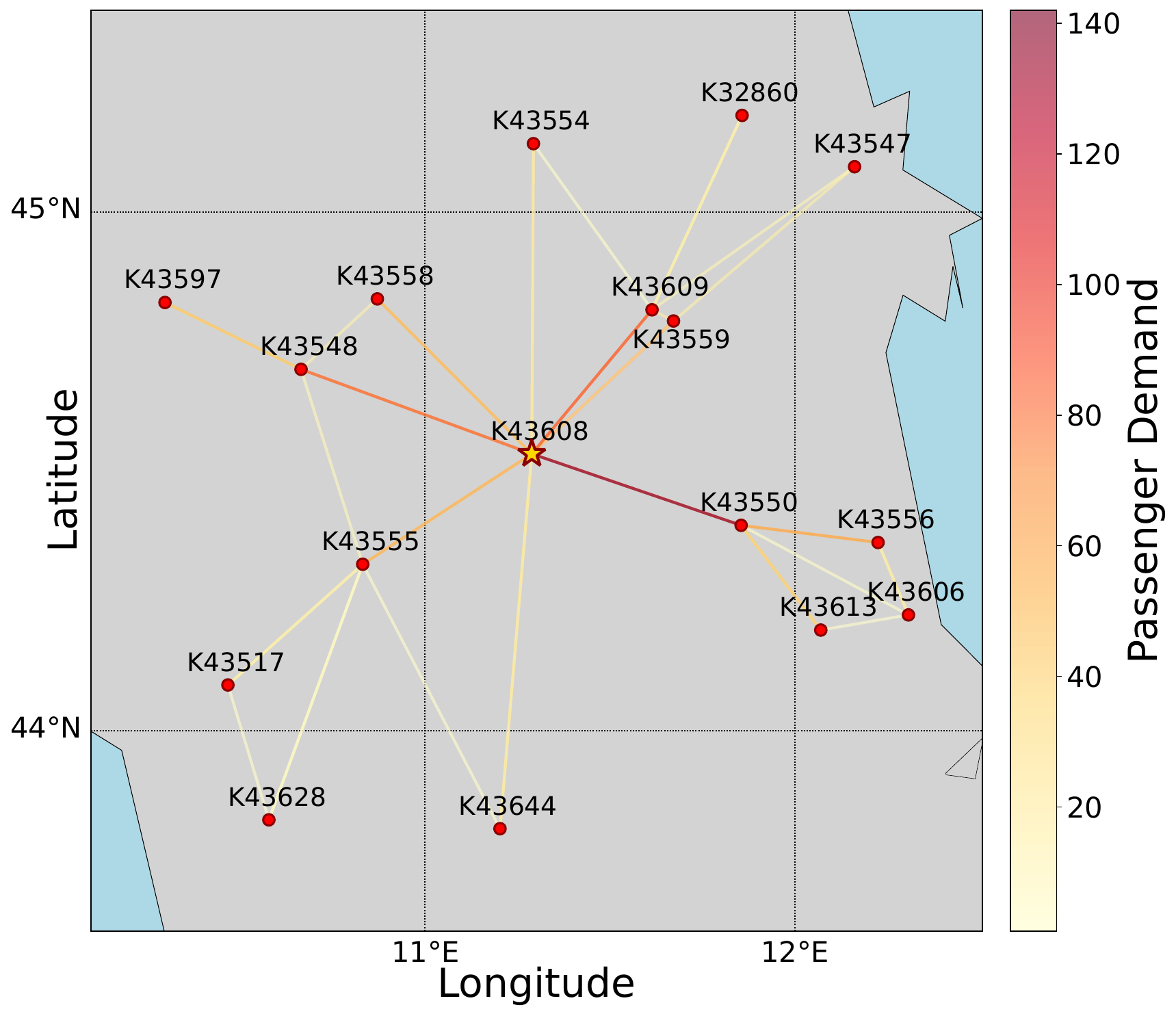}
  \caption{Passenger demand per directed flight leg across the network (passengers aggregated per directed leg).}
  \label{fig:demand_map}
\end{figure}

\begin{figure}
  \centering
  \includegraphics[width=\linewidth]{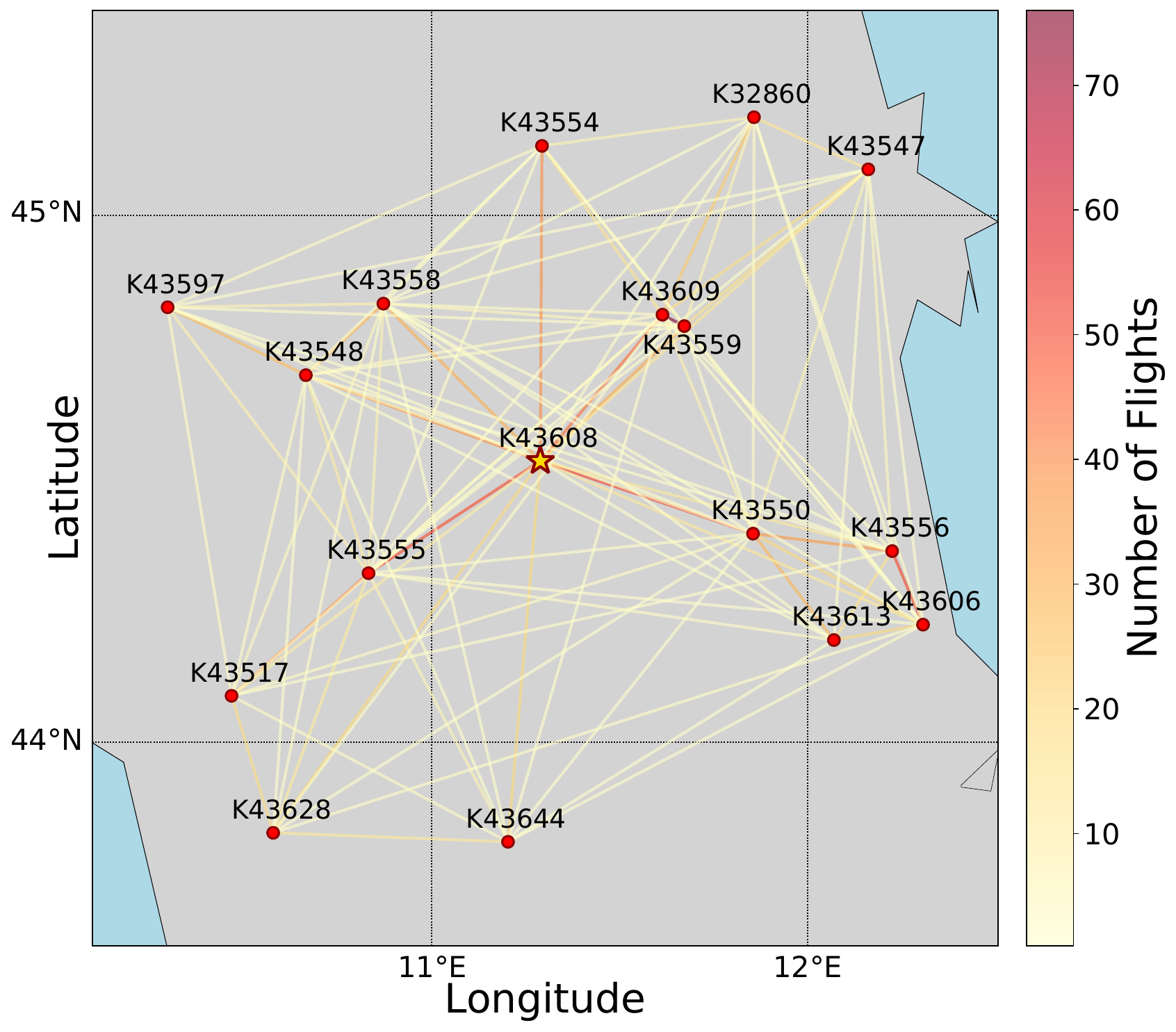}
  \caption{\gls{iam} service frequency per leg derived from the processed \gls{iam} schedule (number of \gls{iam} operations per directed leg).}
  \label{fig:iam_map}
\end{figure}

As shown in \cref{fig:demand_map}, the passenger--demand map consolidating counts per directed leg reveals a pronounced core around the 44--45$^\circ$N and 11--12$^\circ$E band, with the densest cluster of labeled nodes—including K43608 and nearby nodes such as K43555, K43548, K43609, and K43550—supporting the highest intensities. The legend indicates a broad range of volumes (reaching up to $\sim$140 passengers on the busiest legs), signaling several high-demand corridors within this central area, while peripheral nodes (e.g., K32860, K43644) delineate the edge of the demand footprint. Motivated by this spatial concentration of demand, we next examine whether the \gls{iam} supply is co-located with these corridors by turning to the service-frequency map in \cref{fig:iam_map}.

The \gls{iam} network map in \cref{fig:iam_map} summarizes the frequency of operations per leg after processing the \gls{iam} schedule, highlighting a core of repeated operations centered on K43608 and adjacent nodes. The legend shows that some connections accumulate up to $\sim$70 flights over the simulation horizon, producing a markedly denser mesh of \gls{iam} activity in the central cluster, with frequencies tapering toward the periphery (e.g., K32860, K43644). Building on this supply view, we now compare \cref{fig:demand_map} and \cref{fig:iam_map} to assess demand--supply alignment across the network.

A side-by-side reading of \cref{fig:demand_map} and \cref{fig:iam_map} indicates a broad alignment between where passenger demand is strongest and where \gls{iam} service is most frequent—particularly within the central cluster anchored by K43608. At the same time, visual inspection suggests pockets of potential mismatch near the network’s edges, where demand is present but service frequency appears comparatively sparse, and conversely, legs that receive repeated operations despite modest demand. These spatial patterns point to clear operational levers: computing a per-leg balance metric (e.g., passengers per \gls{iam} flight) and ranking legs by under-/over-service would identify candidates for retiming or capacity adjustments, while a demand--vs.--frequency scatterplot could quantify the correlation and isolate outliers for targeted intervention.

\ParamBarFigure
  {data/speed.csv}
  {speed}
  {satisfaction}
  {Speed (m/s)}
  {Demand-met rate vs.\ speed}
  {fig:satisfaction_by_speed}

\ParamBarFigure
  {data/staying.csv}
  {staying}
  {satisfaction}
  {Staying (fraction)}
  {Demand-met rate vs.\ staying}
  {fig:satisfaction_by_staying}

To assess the robustness of the proposed \gls{iam} feeder system, we vary two key parameters: the \gls{iam} cruise speed and the passenger diversion (“staying”) parameter, which specifies the share of arriving passengers at \gls{blq} who may use \gls{iam} instead of connecting to another flight. \gls{iam} fleet size and the vehicle capacity have been studied in our earlier work \cite{gui2025regional}. All other settings remain fixed to the baseline scenario: demand is derived from the OAG schedule of 5 May 2024 (scaled with a 0.8 load factor), evaluated in 5‑minute rolling windows between 09:00 and 17:00. The \gls{iam} fleet consists of 50 vehicles with a capacity of 10 seats, operating with a default speed of 30 m/s and a 2‑hour pre‑departure buffer for feeder passengers. The baseline diversion rate is set to 70\%, and only this parameter and the \gls{iam} speed are varied in the sensitivity study in order to isolate their impact on coverage, timing feasibility, and overall \gls{iam} performance, measured through the demand‑met rate.

We first examine the effect of \gls{iam} cruise speed, as illustrated in \cref{fig:satisfaction_by_speed}. Here, the system exhibits very strong responsiveness: increasing the operational speed from 10 m/s to about 25–30 m/s generates a steep improvement in the demand‑met rate (from roughly 40\% to nearly 90\%). Beyond 35 m/s, the curve plateaus and eventually reaches values above 99\%. This steep early increase reflects that flight time dominates the end‑to‑end service duration at low speeds, and that even modest accelerations greatly enhance fleet productivity and spatial reach. The plateau at higher speeds indicates that once flight time is no longer the limiting factor, other components of the system—such as vehicle availability or node‑level congestion—become the effective bottlenecks. Consequently, while speed is the most influential parameter among the three, its effect saturates once the system operates above a critical velocity threshold.

Finally, the staying parameter represents the share of arriving passengers at \gls{blq} who immediately exit the airport system to reach nearby cities in Emilia‑Romagna rather than continuing through the broader network. As this proportion increases, a larger fraction of requests concentrates on short \gls{blq}→nearby‑city legs. Given a fixed \gls{iam} fleet, that concentration raises local utilization and queuing on these short spokes, effectively tightening the capacity bottleneck around the hub. Consistent with this mechanism, the demand‑met rate declines slightly as the staying parameter rises—from about 0.5 to 0.9 on the x‑axis, the served share drops by roughly one to two percentage points—indicating that more “staying” traffic diverts constrained vehicle time to a subset of legs and leaves marginal requests unserved. In other words, higher staying proportions intensify local competition for limited vehicles and modestly erode overall satisfaction, as observed in \cref{fig:satisfaction_by_staying}.

Taken together, these two analyses show a coherent hierarchy of sensitivity: \gls{iam} speed is the dominant factor, enabling large performance jumps; the staying proportion exerts only marginal effects within the tested range. These insights inform which parameters merit prioritization during system calibration and which offer limited operational leverage.

The sensitivity analysis indicates that the main conclusions are resilient to temporal perturbations. Changes in \gls{iam} cruise speed and passenger diversion mainly affect the share of demand that can be served within required time windows, rather than the qualitative feasibility of \gls{iam} feeder operations. Both schedule uncertainty and weather impacts can be interpreted as effective reductions in available temporal margins or achievable speeds, which are naturally handled by the rolling‑window structure. Accordingly, the observed sensitivity to \gls{iam} speed provides a first‑order proxy for robustness under delayed or weather‑constrained operations, while the structural benefits of \gls{iam} for hub accessibility remain unchanged.

\section{Conclusion} \label{sec:conclusion}

This paper examined Innovative Air Mobility as a dedicated hub feeder and dispersal service for \gls{blq}, using an \gls{iam}‑only time‑expanded network in which \gls{blq}’s flight schedule is employed solely to time‑stamp demand. Unlike multimodal \gls{iam} studies, our formulation removes ground modes and commercial flight arcs from the network, thereby isolating \gls{iam}’s contribution to hub accessibility. The framework couples passenger routing (A*) and fleet dispatching (dynamic programming) in a sequential, demand‑consistent pipeline, and introduces a rolling‑window demand model that aligns \gls{iam} requests with \gls{blq} arrivals and departures via post‑arrival and pre‑departure buffers.

Using real OAG schedules, we quantified \gls{iam}’s potential to reduce access times between \gls{blq} and nearby airports, documented substantial time advantages over rail for multiple links, and synchronized \gls{iam} operations with hub‑driven peaks. We also corrected for partial schedule coverage with a daily scaling factor and analyzed sensitivity to \gls{iam} speed and passenger diversion, showing how these parameters shape coverage and temporal feasibility while other settings remain fixed.

While the proposed framework relies on scheduled operations, real‑world airport environments are subject to stochastic disturbances such as flight delays, weather disruptions, and variability in vertiport turnaround times. These factors would primarily affect the temporal alignment between flight‑driven demand and \gls{iam} availability rather than the underlying structure of the \gls{iam} network. In the present study, demand is processed through short rolling time windows, which naturally limits the propagation of localized schedule perturbations and allows moderate delays to be absorbed through window shifts. Moreover, the substantial travel‑time margins observed between \gls{iam} and rail for most airport pairs provide inherent slack that can mitigate the impact of timing uncertainty. Although explicit stochastic modeling is beyond the scope of this work, the sequential nature of the passenger routing and dispatch mechanisms would allow the framework to be extended toward real‑time or receding‑horizon re‑optimization as updated information becomes available. 

Future work will enrich the operational realism of the demand model (e.g., differentiated diversion by market segment), benchmark the heuristic pipeline against exact MILP instances on reduced networks, and extend the analysis to weather, airspace, and infrastructure constraints. Overall, the results support the viability of \gls{iam} as a short‑range connector that strengthens \gls{blq}’s regional reach without modifying airline schedules.

\section*{Acknowledgment}

This research was co-funded by SESAR Joint Undertaking and European Commission, within the project “Implemented Synergies, data sharing contracts and goals between transport modes and air transportation” (SIGN-AIR) under grant number 101114845, and the project “Planning Regional-Scale Multimodal Operations for Innovative Air Mobility Services” (PRIAM) under grant number 101167262. 

\bibliographystyle{IEEEtran}
\bibliography{mybibfile}

\end{document}

%% file: tex/param-barplots.tex
% ...

\newif\ifbarplottighty
\barplottightytrue

\newcommand{\ParamBarFigure}[6]{%
  \pgfplotstableread[col sep=comma]{#1}\datatable
  \pgfplotstablesort[sort key={#2}]\sortedtable{\datatable}

  \begin{figure}[t]
    \centering
    \begin{tikzpicture}
      \begin{axis}[
        width=0.92\linewidth,
        height=0.50\linewidth,
        ybar,
        bar width=10pt,
        enlarge x limits=0.08,
        enlarge y limits={upper, value=0.15},
        xlabel={#4},
        ylabel={Demand-met rate},
        ymin=0, ymax=1.02,
        scaled y ticks=false,
        yticklabel={
          \pgfmathparse{100*\tick}%
          \pgfmathprintnumber[fixed,precision=0]{\pgfmathresult}\%
        },
        ymajorgrids=true,
        grid style={dashed,gray!40},
        tick style={black},
        axis line style={black},
        label style={font=\small},
        tick label style={font=\small},
        title style={font=\small},
        every axis plot/.append style={
          fill=blue!65!white,
          draw=black,
          line width=0.4pt
        },
        nodes near coords,
        nodes near coords align={vertical},
        nodes near coords style={font=\scriptsize},
        point meta=y,
        nodes near coords={
          \pgfmathprintnumber[fixed,precision=2]\pgfplotspointmeta
        },
        xtick=data,
        xticklabel style={
          /pgf/number format/.cd,
          fixed,
          precision=0,
          1000 sep={}
        },
      ]

        \addplot+[ybar] table[
          x expr=\coordindex,
          y={#3}
        ] {\sortedtable};

        \pgfplotstablegetrowsof{\sortedtable}
        \pgfmathtruncatemacro{\NROWS}{\pgfplotsretval-1}
        \pgfplotsset{xtick={0,...,\NROWS}}
        \pgfplotsset{xticklabels from table={\sortedtable}{#2}}

      \end{axis}
    \end{tikzpicture}
    \caption{#5}
    \label{#6}
  \end{figure}
}